%% file: main.tex
\documentclass[11pt]{amsart}
\usepackage{amssymb}
\usepackage{fullpage}

\usepackage{mathtext}
\usepackage[OT2,OT1]{fontenc}

\def\cyr{\fontencoding{OT2}\fontfamily{wncyr}\selectfont}

\def\Ch{\textrm{\cyr CH}}

\input{macros.tex}
\begin{document}

\bigskip\bigskip

\title[Examples of Feigenbaum Julia sets with small Hausdorff dimension]
{Examples of  Feigenbaum Julia sets\\ with small Hausdorff dimension}
\author {Artur Avila and Mikhail Lyubich}

\address{
Laboratoire de Probabilit\'es et Mod\`eles al\'eatoires\\
Universit\'e Pierre et Marie Curie--Boite courrier 188\\
75252--Paris Cedex 05, France
}
\email{artur@ccr.jussieu.fr}

\address{
Mathematics Department and IMS\\
SUNY Stony Brook\\
Stony Brook, NY 11794, USA
}
\email{mlyubich@math.sunysb.edu}

\address{
Department of Mathematics\\
University of Toronto\\
Ontario, Canada M5S 3G3
}
\email{misha@math.toronto.edu}


\dedicatory{To Bodil Branner on her 60th birthday}
\date{\today}

\begin{abstract}
We give examples of infinitely renormalizable quadratic polynomials $F_c: z\mapsto z^2+c$
with stationary combinatorics whose Julia sets have Hausdorff dimension arbitrary close to 1. 
The combinatorics of the renormalization involved  is close to the Chebyshev one. 
The argument is based upon a new tool, 
a ``Recursive Quadratic Estimate'' for the Poincar\'e series of an infinitely renormalizable map. 
\end{abstract}

\setcounter{tocdepth}{1}

\maketitle
\thispagestyle{empty} \input{imsmark}
\SBIMSMark{2004/04}{July 2004}{}

 \input{intro.tex}
\input{chebnew}

\input{cheb.tex}

\input{bib}
\end{document}

%% file: macros.tex
\newcommand{\crit}{{\mathrm{cr}}}

\let\newpf\proof \let\proof\relax 
\newenvironment{pf}{\newpf[\proofname]}{\qed\endtrivlist}

\def\area{\operatorname {area}}

\def\be{\begin{equation}}
\def\ee{\end{equation}}

\def\ba{{\begin{align}}}
\def\ea{{\end{align}}}

\def\u{{\mathbb U}}

\def\op{\overline\partial}

\def\0{{\mathbf 0}}

\def\cal{\mathcal}

\newtheorem{thm}{Theorem}[section]
\newtheorem{cor}[thm]{Corollary}

\newtheorem{lemma}[thm]{Lemma}

\theoremstyle{remark}
\newtheorem{rem}{Remark}[section]

\numberwithin{equation}{section}

\def \bn {\hfill \\ \smallskip\noindent}

\theoremstyle{definition}

\def\proof{\bn {\bf Proof.} }

\def\note#1
{\marginpar
{\tiny $\leftarrow$
\par
\hfuzz=20pt \hbadness=9000 \hyphenpenalty=-100 \exhyphenpenalty=-100
\pretolerance=-1 \tolerance=9999 \doublehyphendemerits=-100000
\finalhyphendemerits=-100000 \baselineskip=6pt
#1}\hfuzz=1pt}

\newcommand{\ra}{\rightarrow}

\def\msk{\medskip}

\def\sm{\smallsetminus}

\def\ssm{\smallsetminus}
\renewcommand{\setminus}{\ssm}

\newcommand{\diam}{\operatorname{diam}}

\newcommand{\dist}{\operatorname{dist}}

\renewcommand{\mod}{\operatorname{mod}}

\newcommand{\HD}{\operatorname{HD}}

\newcommand{\meas}{\operatorname{meas}}

\newcommand{\Dil}{\operatorname{Dil}}

\newcommand{\eps}{{\epsilon}}

\newcommand{\de}{{\delta}}
\newcommand{\la}{{\lambda}}

\newcommand{\FF}{{\cal F}}

\newcommand{\OO}{{\cal O}}

\newcommand{\C}{{\mathbb C}}
\newcommand{\D}{{\mathbb D}}

\newcommand{\R}{{\mathbb R}}

\newcommand{\V}{{\mathbb V}}
\renewcommand{\U}{{\mathbb U}}

\def\B0{{\bold{0}}}

\catcode`\@=12

\def\Empty{}
\newcommand\oplabel[1]{
  \def\OpArg{#1} \ifx \OpArg\Empty {} \else
  	\label{#1}
  \fi}
		
\newcommand{\comm}[1]{}
\newcommand{\comment}[1]{}

%% file: imsmark.tex
\def\IMSmarkvadjust{0 pt}
\def\IMSmarkhadjust{0 pt}
\def\IMSmarkhpadding{0 pt}
\def\IMSpubltext{Published in modified form:}
\def\SBIMSMark#1#2#3{
 \font\SBF=cmss10 at 10 true pt
 \font\SBI=cmssi10 at 10 true pt
 \setbox0=\hbox{\SBF \hbox to \IMSmarkhpadding{\relax}
                Stony Brook IMS Preprint \##1}
 \setbox2=\hbox to \wd0{\hfil \SBI #2}
 \setbox4=\hbox to \wd0{\hfil \SBI #3}
 \setbox6=\hbox to \wd0{\hss
             \vbox{\hsize=\wd0 \parskip=0pt \baselineskip=10 true pt
                   \copy0 \break%
                   \copy2 \break%
                   \copy4 \break}}
 \dimen0=\ht6   \advance\dimen0 by \vsize \advance\dimen0 by 8 true pt
                \advance\dimen0 by -\pagetotal
	        \advance\dimen0 by \IMSmarkvadjust
 \dimen2=\hsize \advance\dimen2 by .25 true in
	        \advance\dimen2 by \IMSmarkhadjust

%
%
  \openin2=publishd.tex
  \ifeof2\setbox0=\hbox to 0pt{}
  \else 
     \setbox0=\hbox to 3.1 true in{
                \vbox to \ht6{\hsize=3 true in \parskip=0pt  \noindent  
                {\SBI \IMSpubltext}\hfil\break
                \input publishd.tex 
                \vfill}}
  \fi
  \closein2
  \ht0=0pt \dp0=0pt
 \ht6=0pt \dp6=0pt
 \setbox8=\vbox to \dimen0{\vfill \hbox to \dimen2{\copy0 \hss \copy6}}
 \ht8=0pt \dp8=0pt \wd8=0pt
 \copy8
 \message{*** Stony Brook IMS Preprint #1, #2. #3 ***}
}

%% file: intro.tex
\section{Introduction}

One of the most remarkable objects in complex dynamics are the fixed points
of the Douady-Hubbard renormalization operator.  Such objects have a
distinguished place in the dictionary between rational maps and Kleinian
groups (see \cite {McM2}).
Existence of the renormalization fixed points
established in the works of Sullivan \cite{S} and McMullen \cite{McM2} 
(under certain  assumptions)
implies many beautiful features (self-similarity, universality, hairyness,...) 
of Feigenbaum Julia sets 
(see \S \ref{defs} for the definition).
However, even with this thorough information, some basic questions concerning
measure and dimension of these Julia sets have remained unsettled. 

One of the key questions (asked, for instance, in \cite
{McM2}) regarding the geometry of Feigenbaum Julia sets has been the following: 
{\it Is the Hausdorff dimension of a Feigenbaum Julia set always equal to $2$}?
In \cite {AL} we supply a fairly large class of Feigenbaum Julia sets with $\HD(J)<2$,
thus giving a negative answer to the above question.
In this paper we show that in fact the dimension of a Feigenbaum Julia set 
can be arbitrarily close to 1:
\smallskip

\begin{thm}\label{HD}
 There exists a sequence of  Feigenbaum quadratic polynomial $F_{c_p}: z\mapsto z^2 + c_p$
with $c_p\in \R$, $c_p\to -2$, 
such that  $\HD(J(F_{c_p}))\to 1$ as $p\to \infty$. 
\end{thm}

\smallskip

Hausdorff dimension is closely related to another geometric characteristic of the Julia set,
the {\it critical exponent} $\de_\crit$ (see \S \ref{P-series}).
In fact, for a Feigenbaum map $F_c$, 

$$
   \HD(J(F_c))=\de_\crit(J(F_c)),
$$

\noindent 
provided $\meas(J(F_c))=0$  \cite{AL},
and the same is true for the associated renormalization fixed point $f_c$.%
\footnote{In fact, it follows from Bishop's work  \cite{B} that for any $c$, 
 $ \HD(J(F_c))\,\leq \,\de_\crit(J(F_c))$, provided $\meas(J(F_c))=0$.}     
This allows us to reduce Theorem \ref{HD} to the following two results.

\smallskip

\begin{thm} \label {thmA}

Let $f_p$ be the fixed point of the renormalization operator of period $p$
with combinatorics closest to the Chebyshev one.  Then
$\delta_\crit(f_p)\to 1$ as $p \to \infty$. 

\end{thm}

\smallskip

The proof of this theorem is based upon a ``Recursive Quadratic Estimate'' for the Poincar\'e series
which provides a new efficient tool for getting bounds on the critical exponent. 

\smallskip

\begin{thm} \label {thmB}

For large $p$, $\area J(f_p)=0$.

\end{thm}

\msk
{\it Remarks.} (1)
The class of Feigenbaum maps with $\HD(J(f))<2$ 
supplied in \cite{AL} is qualitatively the same as the class treated in 
Yarrington's thesis \cite{Y} (see also \S 9 of \cite {AL})  for which $\area(J)=0$
(which in turn, is qualitatively the same as the class of infinitely renormalizable maps for which
{\it a priori} bounds were established in \cite{puzzle}).  
Though Theorem \ref{thmB} is not formally covered by \cite{AL,Y},
it is proved by a similar method, which  becomes more direct in our situation. 
Similarly,  to prove Theorem \ref{thmA} we adjust the method of \cite{AL} 
to the Chebyshev combinatorics, which makes it (in this combinatorial case)
simpler and more powerful.

\msk
(2)
There are many quadratic polynomials whose Julia sets have Hausdorff dimension two \cite{Sh}. 
However, it is still unknown whether there exist Feigenbaum Julia sets with this property.

\comm{
\begin{pf}

To conclude this corollary from Theorems~\ref {thmA} and~\ref {thmB},
we can use either the previously mentioned Theorem~C of \cite {AL}) or
the result of Bishop \cite {B}, which shows that $\area J(f)=0$ implies the
inequality $\HD(J(f)) \leq \delta_\crit(f)$
for a very general class of maps (including Feigenbaum maps).
\end{pf}
}

\bigskip

{\bf Acknowledgment.}
This work was partly supported by the the Clay Mathematics Institute,
the Guggenheim Foundation, NSF and NSERC. 
It was carried out during authors' meetings at the IMS at Stony Brook, at the  University of Toronto,
and the Institut Henri Poincar\'e. The authors thank all these Institutions and Foundations
for hospitality and support.

\bigskip

%% file: chebnew.tex
\section{Basic concepts}

\subsection{Notations}
$\D_r(z) \equiv \{w \in \C,\, |w-z|<r\}$, $\D_r \equiv \D_r(0)$.
A domain is a connected open subset of $\C$.  A topological disk is a simply
connected domain.
$U \Subset V$ means that $U$ is {\it compactly contained} in $V$.

Notation $a \asymp b$ means that $C^{-1}<a/b<C$ with a constant $C>0$ independent of
particular $a$ and $b$ under consideration; 
$a \approx b$ means that $a$ is close to $b$.

We usually denote the $p$-fold iterate of a map $f$ by $f^p$,
but occasionally use a more forceful notation $f^{\circ p}$.

Let $\omega(x) \equiv \omega_f(x)=\cap_{m \geq 0} \overline {\{f^k(x),\, k \geq m\}}$ 
denote the {\it $\omega$-limit set} of $x$.

For a quadratic-like map $f: U\ra V$ (see below) with the critical point at 0,
let  $\OO(f) \equiv \overline {\{f^k(0),\, k>0\}}$
denote its {\it postcritical set}.


\smallskip

\subsection{Quadratic-like maps and renormalization}\label{defs}

\smallskip

A {\it quadratic-like map} is a holomorphic double covering map
$f:U \to V$ where $U,V \subset \C$ are topological disks and
$U \Subset V$.  Such a map has a unique critical point which we
will assume to be $0$.  Let $K(f) \equiv \cap_{k=0}^\infty f^{-k}(U)$
denote the {\it filled Julia set} of
$f$ and let $J(f) \equiv \partial K(f)$ denote its {\it Julia set}.

Two quadratic-like germs $f$ and $g$ are said to be {\it hybrid equivalent}
if there exists a quasiconformal map $h:\C \to \C$ satisfying $h(f(x))=g(h(x))$ 
for $x$ near $J(f)$ such that $\op h|J(f)=0$. Any quadratic-like map $f: U\ra V$ with connected 
Julia set is hybrid equivalent to a unique quadratic polynomial $F: z\mapsto z^2+c$ 
called the {\it straightening} of $f$ \cite{DH}. Moreover, the dilatation $\Dil(h)$ of the 
(appropriately chosen) conjugacy $h$ depends only on $\mod(V\sm U)$, and $\Dil(h)\to 1$ as
 $\mod(V\sm U)\to \infty$.

\bigskip
The Julia set $J(f)$ of a quadratic-like map is either connected or Cantor.
If $J(f)$ is connected, there exists a unique repelling or parabolic 
fixed point $\beta= \beta(f) \in J(f)$ such that $J(f) \setminus \{\beta(f)\}$ is connected. 
The other fixed point is denoted by $\alpha=\alpha(f)$.  
We will only consider quadratic-like maps with connected Julia set.

A quadratic-like map which is considered only up to choice of domains is
called a {\it quadratic-like germ}.  More precisely, one says that two
quadratic-like maps with connected Julia sets represent the same germ if
they have a common Julia set and coincide in a neighborhood of it.  We shall
consider quadratic-like germs up to affine conjugacy.

\bigskip

A quadratic-like map $f: U\ra V$ is called {\it renormalizable} with period $p>1$ if there
exist topological disks $U'\Subset V'$ containing the critical point such that
\begin{enumerate}
\item 
 $g\equiv f^p: U'\ra V'$ is a quadratic-like map with connected Julia set $J(f')$ called 
 a {\it pre-renormalization} of $f$;
\item
 For every $1 \leq k \leq p-1$, either $f^k(J(g)) \cap J(g)=\emptyset$ 
 or $f^k(J(g)) \cap J(g)=\{\beta(g)\}$.
\end{enumerate}
The renormalization operator $R$ is defined on the space of germs by letting $Rf=g$. 
The minimal $p=p(f) > 1$ for which $f$ is renormalizable 
is called the {\it renormalization period} of $f$.
In what follows, the operator $R$ will always correspond the this minimal period. 

A quadratic-like map $f:U \to V$ is said to be a {\it renormalization fixed point}
if $f$ is renormalizable and $Rf=f$.  In other words, 
$f^p(x)=\lambda f(\lambda^{-1} x)$ near $J(g)$ for some 
$\lambda \in \D \setminus \{0\}$,
where $p$ is the renormalization period of $f$ and $g$ is a
pre-renormalization of $f$.

\subsection{Poincar\'e series}\label{P-series}

Let $f:U \to V$ be a quadratic-like map.

Sullivan's Poincar\'e series  \cite{S-conformal} is defined as follows: 
$$
\Xi_\delta(z)=\sum_{k=0}^\infty \sum_{f^k(w)=z} |Df^k(w)|^{-\delta},\quad
     z \in V \setminus \OO(f),\quad \delta>0.
$$
It follows from the Koebe Distortion Theorem that $\Xi_\delta(z) \leq
C(z,z')^\delta \Xi_\delta (z')$ for any $z,z' \in V \setminus \OO(f)$.
In particular, $\Xi_\delta$ is finite or infinite independently  of $z$.

The function $\delta \mapsto \Xi_\delta$ is obviously convex.
By definition, the {\it critical exponent}, $\delta_\crit(f) \in [0,\infty]$, 
is the unique value of $\de$  that separates convergent $\Xi_\delta$ from  divergent ones.
The critical exponent depends only on the germ of $f$ near $K(f)$.

It is easy to see that $\Xi_2$ is always finite (area argument) and,
since $J(f)$ is assumed to be connected,
$\Xi_1=\infty$ (length argument), see \S 2.9 of \cite {AL}.
Thus we actually have $\delta_\crit(f) \in [1,2]$.%
\footnote{In fact, $\delta_\crit>1$, unless $J(f)$ is a real analytic curve.}  

\smallskip

\subsubsection{Poincar\'e series for subfamilies of orbits}

\smallskip

An {\it orbit of length} $k \geq 0$ is a sequence $(x_0,...,x_k)$, 
where $x_k \in V$ and $f(x_i)=x_{i+1}$ for $0 \leq i <k$.
An orbit of zero length is called {\it trivial}.

Given a family $\FF$ of orbits $(x_0,...,x_k)$,
we define a function $\C \to [0,\infty]$
$$
\Xi_\delta(\FF)(z)=\sum_{k=0}^\infty
\sum_{(x_0,...,x_k=z) \in \FF} |Df^k(x_0)|^{-\delta}
$$
(to keep notation shorter, we do not explicitly mention  $f$). 
 Let  $\Xi^{[j]}_\delta$ denote the truncation of $\Xi_\delta$ at level $j$,
$$
\Xi^{[j]}_\delta(\FF)(z)=\sum_{k=0}^j
\sum_{(x_0,...,x_k=z) \in \FF} |Df^k(x_0)|^{-\delta},
$$
with convention that $\Xi^{[j]}=0$ for $j<0$.
Note that $\Xi^{[0]}(\FF)$ is equal to $1$ or $0$ 
depending on whether $\FF$ contains the trivial orbit or not.

\subsubsection{Arrow notation}

Let us introduce a convenient notation for certain families of orbits.
Let $D,E\subset V$, $ S\subset U $.   
By $D \xleftarrow{} E$, we will understand the family of orbits
$(x_0,...,x_k)$ with $x_0 \in E$ and $x_k \in D$.
The family of orbits $(x_0,...,x_k)$ with $x_0
\in E$, $x_k \in D$ and $x_1,...,x_{k-1} \in S$ will be denoted $D
\xleftarrow[F]{} E$.  A ``plus sign'' over the arrow will indicate that
only non-trivial orbits are considered.
The juxtaposition of arrows will denote composition in the
natural way.  For instance, 
$$
D \xleftarrow[S]{+} D \xleftarrow[S]{} E,
$$
denotes the family of orbits $(x_0,...,x_k)$, with $x_0 \in E$, $x_k \in D$,
such that  $x_i \in D$ for some  $0 \leq i<k$, 
and $x_1,...,x_{i-1},x_{i+1},...,x_{k-1} \in S$.

\section{Quadratic recursive estimate}

We will now introduce a version of the Quadratic Recursive Estimate 
which is sufficient for purposes of this paper  
(see \cite {AL} for a finer version). 
We shall restrict ourselves to the case of renormalization fixed points.
The argument is based on a combinatorial breakdown of  orbits
which exploits the scaling self-similarity of the dynamics.

Let $f:U \to V$ be a fixed point of renormalization
of period $p$, $f^p(x)=\lambda f(\lambda^{-1}x)$ near $0$.  
Let $U'=\lambda U$, $V'=\lambda V$, and  let $gx \equiv f^p:U' \to V'$.
Let $A=V \setminus U$, $A'=V' \setminus U'$.
We assume that $V' \subset U$,  $g$ is the first return from $U'$ to $V'$,
and that $\OO(f)$ does not intersect $\bar A'$. 

\begin{lemma} \label {recursion}

Let $s_j(\delta)=\sup_{z\in A'} \, \Xi_\delta^{[j]}(A' \xleftarrow{} U)(z)$.
Then
$$
s_{j+1}(\delta) \leq P_\delta(s_j(\delta)),
$$
where $x \mapsto P_\delta(x)$ is a quadratic polynomial with
positive coefficients which can be  expressed explicitly in terms of the Poincar\'e series
$\Xi_\de(\FF)$ over families $\FF$ of orbits that do not accumulate on 0. 

If $P_\delta$ has a positive fixed point $s$
then 
$$
\sup_{z\in A'}\, \Xi_\delta(A' \xleftarrow{} U)(z)=\lim s_j \leq s,
$$
so that  $\delta_\crit(f) \leq \delta$. 
\end{lemma}

\begin{pf}
In what follows, the $\sup$ is always taken over $z$, 
the terminal point of the orbit in question.
We will also omit the truncation parameter ($j$ or $j+1$) in the notation.

We can decompose $A \xleftarrow{} U$ into two groups:
$A \xleftarrow[U \setminus V']{} U \setminus V'$ and
$A \xleftarrow[U \setminus V']{} A' \xleftarrow{} U$.
This gives the inequality
\be \label {firstine}
\sup\, \Xi_\delta(A \xleftarrow{} U) \leq \sup\, \Xi_\delta(A \xleftarrow[U
\setminus V']{} U \setminus V')+\sup\, \Xi_\delta(A \xleftarrow[U \setminus
V']{}
A') \sup\, \Xi_\delta(A' \xleftarrow{} U).
\ee
In turn, we can decompose $A' \xleftarrow{+} U$ into two groups: 
\begin{enumerate}
\item $A' \xleftarrow[U \setminus A']{} U \setminus A'$, which can be
further decomposed into
$$
A' \xleftarrow[U \setminus V']{} U \setminus V', \quad A' \xleftarrow[U
\setminus A']{} U', \quad \text {and} \quad
A' \xleftarrow[U \setminus A']{} U' \xleftarrow[U \setminus V']{}
U' \setminus V'; 
$$
\item $A' \xleftarrow[U \setminus A']{+} A' \xleftarrow{} U$,
which can be further decomposed into
$$
A' \xleftarrow[U \setminus V']{+} A' \xleftarrow{} U \quad \text {and} \quad
A' \xleftarrow[U \setminus A']{} U' \xleftarrow[U \setminus V']{} A'
\xleftarrow{} U. 
$$
\end{enumerate}
This gives the following inequality
\begin{align} \label {secondine}
\sup\, \Xi_\delta(A' \xleftarrow{} U) \leq 
1&+ \sup\, \Xi_\delta(A' \xleftarrow[U \setminus V']{} U \setminus V')\\
\nonumber
& + \sup\, \Xi_\delta (A' \xleftarrow[U \setminus A']{} U') 
  \left( 1 +  \sup\, \Xi_\delta(U' \xleftarrow[U \setminus V']{} U \setminus V')\right)  \\
\nonumber
&+\sup\, \Xi_\delta(A' \xleftarrow[U \setminus V']{+} A')
\sup\, \Xi_\delta(A' \xleftarrow{} U)\\
\nonumber
&+\sup\, \Xi_\delta(A' \xleftarrow[U \setminus A']{} U')
\sup\, \Xi_\delta(U' \xleftarrow[U \setminus  V']{} A')
\sup\, \Xi_\delta(A' \xleftarrow{} U),
\end{align}
where the first term, $1$,  accounts for the trivial orbits.

Notice that since $x \mapsto  f^p x$ is the first return
map from $U'$ to $V'$, if $(x_0,...,x_k)$ belongs to
$A \xleftarrow{} U$ then $(\lambda x_0,...,f^{kp}(\lambda x_0))$
belongs to $A' \xleftarrow[U \setminus A']{} U'$.  This
correspondence is readily seen to be a bijection between
$A \xleftarrow{} U$ and $A' \xleftarrow[U \setminus A']{} U'$
preserving the  weights of the Poincar\'e series.  Hence
$$
\Xi_\delta(A \xleftarrow{} U)(x)=\Xi_\delta(A' \xleftarrow[U \setminus A']{}U')(\lambda x)
$$
and
\be \label {equal}
\sup\, \Xi_\delta(A' \xleftarrow[U \setminus A']{} U')=
\sup\, \Xi_\delta(A \xleftarrow{} U).
\ee

Plugging (\ref {firstine}) into  (\ref {equal}), and
then plugging the resulting expression for $\sup\, \Xi_\delta(A' \xleftarrow[U \setminus A']{} U')$
into the 2nd and 4th lines of  (\ref {secondine}),  we obtain
$$
\sup\, \Xi_\delta(A' \xleftarrow{} U) \leq \alpha+\beta
\sup\, \Xi_\delta(A' \xleftarrow{} U)+\gamma
\sup\, \Xi_\delta(A' \xleftarrow{} U)^2,
$$
where 
\be \label {constantpart}
\alpha = 1+\sup\, \Xi_\delta(A' \xleftarrow[U \setminus V']{} U \setminus V')+
\sup\, \Xi_\delta(A \xleftarrow[U \setminus V']{} U \setminus V')
\left (1+\sup\, \Xi_\delta(U' \xleftarrow[U \setminus V']{} U \setminus V')
\right ),
\ee
\begin{align} \label {linearpart}
\beta = \sup\, \Xi_\delta(A' \xleftarrow[U \setminus V']{+} A')&+
\sup\, \Xi_\delta(A \xleftarrow[U \setminus V']{} A')
\left (1+\sup\, \Xi_\delta(U' \xleftarrow[U \setminus V']{} U \setminus V')
\right )\\
\nonumber
&+\sup\, \Xi_\delta(A \xleftarrow[U \setminus V']{} U \setminus V')
\sup\, \Xi_\delta(U' \xleftarrow[U \setminus V']{} A'), 
\end{align}
   and 
\be \label {quadraticpart}
\gamma = \sup\, \Xi_\delta(A \xleftarrow[U \setminus V']{} A')
\sup\, \Xi_\delta(U' \xleftarrow[U \setminus V']{} A').
\ee
This is the desired quadratic recurrence estimate for
$\sup\, \Xi_\delta(A' \xleftarrow{} U)$. The above three formulas give an explicit
expression of the coefficients $\alpha$, $\beta$ and $\gamma$ of $P_\de$ in terms of
Poincar\'e series over families of orbits that do not accumulate on $0$. 

For the last statement, notice that $s_j \leq P^j(s_{-1})=P^j(0) \leq s$
for all $j$.  Thus, for every $z \in A'$ we have $\Xi_\delta(z) \leq
\sup\, \Xi_\delta(A' \xleftarrow{} U)=\lim_{j \to \infty}
s_j \leq s$, which shows that $\delta_\crit(f) \leq \delta$.
\end{pf}

\section{Renormalization with combinatorics closest to Chebyshev}

In this section we will show that the critical exponent of maps
with combinatorics ``close to Chebyshev'' can be arbitrarily close to $1$.
Our bounds on the critical exponent will be based on direct estimates
of the coefficients of the quadratic recursive polynomial
corresponding to a nearly Chebyshev map.

\subsection{Basic properties}

Let $\Ch(x)=2-x^2$ be the Chebyshev polynomial.
Let $f_p$ {\it be the fixed point of the renormalization operator of period $p$,
with (real) combinatorics closest to Chebyshev}: 
$f_p$ is combinatorially characterized among fixed points of renormalization 
of period $p$ by being
(up to affine conjugacy) a real-symmetric quadratic-like germ such that
$f_p(0)>0$ and $f_p^i(0)<0$, $1<i<p$.  The existence of $f_p$ is a
particular case of a result of Sullivan \cite {MS}.

We normalize $f_p$ so that its orientation preserving fixed point is $-2$.
Let $-1<\lambda_p<0$ be the scaling factor of $f_p$.
Then we have near zero:
$$
  g_p:= f_p^{\circ p}(x)=\lambda_p f_p(\lambda_p^{-1} x). 
$$
Notice that $[-2,2] \subset J(f_p)$ and $f_p:[-2,2] \to [-2,2]$ is
a unimodal map. 
Let $\alpha_p>0$ stand for the orientation reversing fixed point of $f_p$.

A basic fact is that all of
the $f_p$ belong to some fixed {\it Epstein class}, that is, there
exists $\epsilon>0$ such that $f_p:[-2,2] \to [-2,2]$ extends to a
real-symmetric double covering onto the slit plane
$\C \setminus (\R \setminus (-2-\epsilon,2+\epsilon))$.
(The natural topology in such an Epstein class makes it a compact space.)  
This is a consequence of the real {\it a priori} bounds, see \cite {MS}.
This yields a number of nice properties of the maps $f_p$.  
The ones that are relevant for us are summarized in the following lemma: 

\begin{lemma} \label {class}

Let $p \geq 3$, $T'=(-\alpha_p, \alpha_p)$, $\V'=\{z: |z|<\alpha_p\}$,
and let $\U'$ be the component of $f_p^{-p}(\V')$ containing $0$.
Let $\U= \la^{-1} \U'$ and $\V= \la^{-1} \V'$.
Then
\begin{enumerate}
\item $f_p$ extends to a double covering onto the slit plane 
   $\C \setminus (\R \setminus T)$;
\item $f_p \to \Ch$ uniformly in $[-2,2]$ (in particular $f_p(0) \to 2$ and
$\alpha_p\to 1$);
\item the maps $g_p:\U' \to \V'$ and  $f_p:\U \to \V$ are quadratic-like
for $p$ sufficiently large; 
\item $\mod(\V \ssm \overline \U) = \mod(\V' \ssm \overline {\U'}) \to
\infty$;
\item $\lambda_p \to 0$.
\end{enumerate}

\end{lemma}


\begin{pf}


Let $S_k\subset [-2,2]$ be the component of $(f_p|[-2,2])^{ -(p-k) }(T')$
containing $f_p^k (0)$, $k=0,1,\dots, p$.
Since the intervals $[-2, -\alpha_p]$ and $[\alpha_p, 2]$ are monotonically mapped by $f_p$ 
onto $[-2, \alpha_p]$, the maps $f_p: S_k \ra S_{k+1}$ are diffeomorphisms for $k=1,2,\dots, p-1$.
This implies the first assertion by rescaling.

The second assertion follows from the compactness of the Epstein class and
the first assertion.

Moreover,  $|S_1|^{-1/p} \approx \dist(S_1, f_p(0))^{-1/p} \approx 4$,
where ``$4$'' is the multiplier of the orientation preserving
fixed point $-2$ of $\Ch$.  Since $f_p$ belongs to the Epstein class,
the component of $f_p^{-(p-1)}(\V')$ containing $f_p(0)$ is contained
in the round  disk
with diameter $S_1$.  Hence  $(\diam \U')^{-1/p} \approx 2$,
which implies assertions  (3) and (4) for $g_p$.
The corresponding assertions for $f_p$ are obtained by rescaling.

Since $$ \la_p = \frac{\diam J(g_p)}{\diam J(f_p)}\leq \frac{1}{4}
\diam \U' ,$$
assertion (5) follows, too.
\end{pf}


\subsection{Estimates for the coefficients}

We will now use the information provided by Lemma~\ref {class} to give
direct estimates on the coefficients of the quadratic recursive estimate.

The following lemma gives control of expansion along the orbits that
stay away from 0: 

\begin{lemma} \label {1}

For every $x \in \C \setminus \{-2,2\}$, there exists $K=K(x)$ with the
following properties:

\begin{enumerate}

\item If $\Ch^m(y)=x$, $m \geq 1$, then $|D\Ch^m(y)| \geq K 2^m$;

\item For any $\epsilon>0$ and  $p \geq p_0(\epsilon)$, if
$x \in \frac {1} {2} \V$ 
and $f^m_p(y)=x$, $m \geq 1$, with
$f^k_p(y) \notin \D_\epsilon$, $0 \leq k \leq m-1$, then
$|Df^m_p(y)| \geq K (2-\epsilon)^m$.

\end{enumerate}

Moreover, $K$ depends only on the distance from $x$ to $\{-2,2\}$ and
goes to infinity as $x$ goes to infinity.

\end{lemma}

\begin{pf}

Consider the map $T:\C \setminus \D \to \C$, $T(z)=-(z+z^{-1})$
semi-conjugating $z\mapsto z^2$ to $\Ch$;  $T(z^2)=\Ch(T(z))$.  If
$x=T(x')$, $y=T(y')$ and
$\Ch^m(y)=x$ with $m \geq 1$,
then $D\Ch^m(y)=DT(x') DT(y)^{-1} 2^m x' y'^{-1}$.
Since $|y'|= |x'|^{1/2^m} \leq \sqrt{|x'|}$, we have:
\be\label{Ch}
|D\Ch^m(y)| \geq \frac {|DT(x')|} {|DT(y')|} |x'|^{1/2} 2^m.
\ee
Since  $|DT(y')| \leq 2$ for all $y'\in \C\sm \D$ and
$|DT(x')|$ is bounded away from zero for $x$ outside a
neighborhood of $\{-2,2\}$, (\ref{Ch}) implies (1).

Since the dynamics of $\Ch$ outside a neighborhood of $0$ is hyperbolic and
hence H\"older stable, the second statement follows easily from Lemma~4.1.
\end{pf}

For $0<\rho<1$, let $V'=V_\rho' = \D_\rho$ and  $U'=U_{\rho,p}'= f^{-p}_p(V')|0$.
It follows from Lemma~\ref {class} that for $p>p_0(\rho)$, the map 
$g_p= f_p^{\circ p}: U' \to V'$ is a quadratic-like pre-renormalization of $f_p: U\ra V$,
where $U= U_{\rho,p}= \lambda_p^{-1} U'$ and $V= V_{\rho,p}=\lambda_p^{-1} V'$. 
{\it In what follows,  $\rho$ and $p$  will be usually suppressed in the notation.}

\bigskip

The following two lemmas give control of expansion along the orbits 
that originate near 0. 

\begin{lemma} \label {3}

For every $0<\rho \leq 1/10$, $0<\kappa \leq 1/10$,
and $p>p_0(\kappa,\rho)$, we have
\be \label {y}
|f(y)-2| \leq |y|^{2-\kappa}
\ee
for any $y \in \D_{e^{-\kappa^{-2}}} \setminus U'$.

\end{lemma}

\begin{pf}
Since the map $f^{p-1}: [f(0), 2] \ra [f^p(0), -2]$
has bounded distortion, 
$$
    |f(0)-2|\asymp  |D f^{p-1} (f(0))|^{-1}.
$$
Let $W=W(p,\rho)$ be the connected component of $f^{-(p-1)}(V')$ containing
$f(0)$.
Similarly, since the map $f^{p-1}: W\ra \D_\rho$ has bounded distortion,
$$
  \dist (f(0), \partial W) \asymp \rho |D f^{p-1} (f(0))|^{-1}.
$$
Hence for some $\eta>0$,        
$$
  \dist (f(0), \partial W) \geq \eta \rho |f(0)-2|.
$$
It follows that for  $y \notin U'$ we have:  
$2 |y|^2 \geq \eta \rho |f(0)-2|$. 
On the other hand,  since $|f(0)-2|\to 0$ as $p\to \infty$,
we have: $\eta \rho>|f(0)-2|^{\kappa/4}$
for $p>p_0(\kappa,\rho)$. Hence
$2 |y|^2 \geq |f(0)-2|^{1+\kappa/4}$.  
It implies by an elementary calculation that
$$
|y|^{2-\kappa} \geq 2 |y|^2+(2 |y|^2)^{(1+\kappa/4)^{-1}} \geq
|f(y)-f(0)|+|f(0)-2| \geq |f(y)-2|,
$$
provided  $0<\kappa \leq 1/10$ and $|y|<e^{-\kappa^{-2}}$
\end{pf} 

\begin{lemma} \label {2}

For every $\epsilon>0$, $0<\rho<\rho_0(\epsilon)$, 
and for any period $p \geq p_0(\epsilon,\rho)$,  
the following property holds. 
Assume that $y \in A'$
and let $m \geq 2$ be the minimal moment such that 
$|f^m(y)+2|>1/10$. Then
$$
|Df^m(y)| \geq (2-\epsilon)^m.
$$

\end{lemma}

\begin{pf}

A simple consideration of the local dynamics near $-2$
shows that
\be \label {obvious}
|Df^{m-1}(f(y))| \asymp |f(y)-2|^{-1}.
\ee
Hence $m \leq K-\log |f_p(y)-2|/\log \eta_p$,
where $\eta_p=|Df(-2)|$.
Since $\eta_p \to 4$ as $p \to \infty$, we have
$$
(2-\epsilon)^m \leq (2-\epsilon)^{K-\frac {\log |f(y)-2|} {\log
\eta_p}} \leq 2^K |f(y)-2|^{\frac {-1+\kappa} {2}}
$$
for $0<\kappa<\kappa(\epsilon)$ and $p>p_0(\epsilon)$.

Set $\rho=e^{-\kappa^{-2}}$.  By Lemma~\ref {3}, if $p>p_0(\rho)$ then
$y \in V' \setminus U'$ implies (\ref {y}).  On the other hand, 
(\ref {obvious}) and (\ref {y}) yields:
$$
|Df^m(y)|=|Df(y)|\, |Df^{m-1}(f(y))|>K^{-1} |y|\,
|f(y)-2|^{-1} \geq
K^{-1} |f(y)-2|^{\frac {1} {2-\kappa}-1}.
$$
Thus, we just have to check
$$
K^{-1} |f(y)-2|^{\frac {1} {2-\kappa}-1} \geq 2^K
|f(y)-2|^{\frac {-1+\kappa} {2}},
$$
that is,
$$
|f(y)-2|^{\frac {\kappa(1-\kappa)}{4-2\kappa}} \leq \frac{1}{K 2^K},
$$
which follows from (\ref{y}) and $|y|<\rho=e^{-\kappa^{-2}}$, provided
$\kappa$ is small enough.
\end{pf}

Note that we have obtained the same lower bound ($\log 2-\eps $) for the
Lyapunov exponents of orbits that stay away from 0
and for those that originate quite near 0. It is because the multiplier of the
postcritical fixed point $-2$ is big ($2^2-\eps$).

\begin{lemma}\label{prevlem}

For every $\epsilon>0$, $0<\rho<\rho_0(\epsilon)$, and
$p>p_0(\epsilon,\rho)$, we have

\begin{enumerate}

\item If $(x_0,...,x_k) \in (A \xleftarrow[U \setminus V']{} U \setminus U'$)
then $|Df^k(x_0)| \geq K (2-\epsilon)^k$,
where $K=K(p,\rho)\to \infty $ as $p\to \infty$;

\item If $(x_0,...,x_k) \in (V' \xleftarrow[U \setminus V']{} U \setminus U'$)
 then $|Df^k(x_0)| \geq K (2-\epsilon)^k$ for some
absolute $K$.

\end{enumerate}

\end{lemma}

\begin{pf}

Let us deal with the first statement.
Notice that since $\mod(\V \setminus \overline \u)$ is large when
$p$ is large, for fixed $\rho$ we also have
$\lim_{p \to \infty} \mod(V \setminus \overline U)=\infty$.
Since $\mod(U \setminus J(f)) \geq \mod(V \setminus \overline U)$,
we see that the distance $M(p,\rho)$ between $\partial U$ and $0$ satisfies
$\lim_{p \to \infty} M(p,\rho)=\infty$.
Since $x_k \in A$, we have $|x_k| \geq M(p,\rho)$.
If $x_0 \notin V'$ then Lemma~\ref {1} shows that
$|Df^k(x_0)| \geq K(p,\rho) (2-\epsilon)^k$, where
$\lim_{p \to \infty} K(p,\rho)=\infty$.  If $x_0 \in A'$,
we let $2 \leq k_0 \leq k$ be minimal with $|f^{k_0}(x_0)+2|>1/10$.
Then by Lemma~\ref {1},
$|Df^{k-k_0}(f^{k_0}(x_0))| \geq K(p,\rho) (2-\epsilon)^{k-k_0}$,
where $\lim_{p \to \infty} K(p,\rho)=\infty$, and by
Lemma~\ref {2}, $|Df^{k_0}(x_0)| \geq (2-\epsilon)^{k_0}$, so $|Df^k(x_0)|
\geq K(p,\rho) (2-\epsilon)^k$, and the first statement follows.

The second statement is analogous.
\end{pf}

\begin{lemma} \label {prevlemm}

Let $\delta>1$.  Then
$$
\lim_{p \to \infty}\,
\sup\, \Xi_\delta(A \xleftarrow[U \setminus V']{} U \setminus U')=0, \quad
0<\rho<\rho_0(\delta),
$$
$$
\lim_{\rho \to 0}\,
\limsup_{p \to \infty}\, \sup\, \Xi_\delta(V' \xleftarrow[U \setminus V']{+} A')=0,
$$
$$
\limsup_{p \to \infty}\,
\sup\, \Xi_\delta(V' \xleftarrow[U \setminus V']{} U \setminus V') \leq
K \equiv K(\delta), \quad 0<\rho<\rho_0(\delta).
$$

\end{lemma}

\begin{pf}

By the first statement of Lemma~\ref {prevlem}, for every
$x \in A$ we have
\begin{align*}
\Xi_\delta(A \xleftarrow[U \setminus V']{} U \setminus U')(x) \equiv
\Xi_\delta(\FF)(x) &\leq \sum_{k
\geq 1} \sum_{(x_0,...,x_k=x) \in \FF} |Df^k(x_0)|^{-\delta}\\
&\leq \sum_{k
\geq 1} 2^k K^{-\delta} (2-\epsilon)^{-\delta k}=K^{-\delta} \sum_{k \geq 1}
\left (\frac {2} {(2-\epsilon)^\delta} \right )^k,
\end{align*}
where $K=K(p,\rho)$, $\epsilon=\epsilon(p,\rho)$
satisfy $\lim_{p \to \infty} K(p,\rho)=\infty$,
$\lim_{\rho \to 0} \lim_{p \to \infty} \epsilon(p,\rho)=0$.
The first estimate follows.

Let $m=m(p,\rho)$ be the minimal return time from $A'$ to $V'$.
Then $\lim_{\rho \to 0} \liminf_{p \to \infty} m(p,\rho)=\infty$.
By the second statement of Lemma~\ref {prevlem},
for every $x \in V'$, we have
\begin{align*}
\Xi_\delta(V' \xleftarrow[U \setminus  V']{+} A')(x) \equiv
\Xi_\delta(\FF)(x) &\leq \sum_{k
\geq m} \sum_{(x_0,...,x_k=x) \in \FF} |Df^k(x_0)|^{-\delta}\\
&\leq \sum_{k
\geq m} 2^k K^{-\delta} (2-\epsilon)^{-\delta k}=K^{-\delta} \sum_{k \geq m}
\left (\frac {2} {(2-\epsilon)^\delta} \right )^k,
\end{align*}
where $K$ is an absolute constant and $\epsilon=\epsilon(p,\rho)$
satisfies $\lim_{\rho \to 0}
\lim_{p \to \infty} \epsilon(p,\rho)=0$.  This gives the
second estimate.

By the second statement of Lemma~\ref{prevlem}, for every $x \in V'$ we have
\begin{align*}
\Xi_\delta(V' \xleftarrow[U \setminus  V']{} U \setminus V')(x) \equiv
\Xi_\delta(\FF)(x) &\leq \sum_{k
\geq 1} \sum_{(x_0,...,x_k=x) \in \FF} |Df^k(x_0)|^{-\delta}\\
&\leq \sum_{k
\geq 1} 2^k K^{-\delta}
(2-\epsilon)^{-\delta k}=K^{-\delta} \sum_{k \geq 1}
\left (\frac {2} {(2-\epsilon)^\delta} \right )^k,
\end{align*}
where $K$ is an absolute constant and $\epsilon=\epsilon(p,\rho)$
satisfies $\lim_{\rho \to 0}
\lim_{p \to \infty} \epsilon(p,\rho)=0$.  This gives the
last estimate.
\end{pf}

\noindent {\it Proof of Theorem \ref {thmA}.}
Since any quadratic-like map with connected Julia set satisfies
$\delta_\crit \geq 1$, we only have to show that for every $\delta>1$ and
for every $p$ sufficiently large, $\delta_\crit(f) \leq \delta$, where
$f=f_{p,\rho}$ is some quadratic-like representative of $f_p$.

Fix some $\delta>1$.  By Lemma~\ref {prevlemm}, we can choose
$\rho>0$ so that
$$
\limsup_{p \to \infty}\, \sup\, \Xi_\delta(V' \xleftarrow[U \setminus
V']{+}
A') \leq \frac {1} {4}.
$$

Let $P_\delta$ be the quadratic polynomial defined in Lemma~\ref
{recursion}.  Notice  obvious inequalities
$$
\max \left \{\Xi_\delta(A' \xleftarrow[U \setminus V']{} U \setminus V'),\
\Xi_\delta(U' \xleftarrow[U \setminus V']{} U \setminus V') \right \} \leq
\Xi_\delta(V' \xleftarrow[U \setminus V']{} U \setminus V'),
$$
$$
\max \left \{\Xi_\delta(A \xleftarrow[U \setminus V']{} U \setminus V'),\
\Xi_\delta(A \xleftarrow[U \setminus V']{} A') \right \} \leq
\Xi_\delta(A \xleftarrow[U \setminus V']{} U \setminus U'),
$$
$$
\max \left \{\Xi_\delta(A' \xleftarrow[U \setminus V']{+} A'),\
\Xi_\delta(U' \xleftarrow[U \setminus V']{} A') \right \} \leq
\Xi_\delta(V' \xleftarrow[U \setminus V']{+} A').
$$
By Lemma~\ref {prevlemm}, when $p$ grows, the constant coefficient (\ref {constantpart}) of
$P_\delta$ stays bounded, the linear coefficient (\ref{linearpart}) becomes smaller than $1/3$,
and the quadratic term (\ref{quadraticpart})  goes to $0$.  In particular, for $p$ large $P_\delta$
takes $[0,2 P_\delta(0)]$ into itself, and hence it has a fixed point.  By
Lemma~\ref {recursion}, $\delta_\crit(f) \leq \delta$ as desired.
\qed

\comm{
\begin{rem}

It is easy to see that we also have
$$
\lim_{p \to \infty} \delta_{\crit}(g_p)=1,
$$
where $g_p$ is the {\it straightening} of $f_p$, that is, the unique
(up to affine conjugacy) quadratic polynomial in the same {\it hybrid
class}\footnote{Recall that two quadratic-like germs $f$ and $g$ are said to
belong to the same hybrid class if there
exists a quasiconformal map $h:\C \to \C$
satisfying $h(f(x))=g(h(x))$ for $x$ near $J(f)$ and such that $\op
h|J(f)=0$.} of $f_p$.

Indeed, the sequence $f_p$ is converge (in the natural topology of
quadratic-like germs) to the Chebyshev polynomial which belongs to the
quadratic family.  Thus, for $p$ large,
$f_p$ and $g_p$ are close quadratic-like germs in the same hybrid class
which easily implies (see for instance Lemma 3.12
\note {We should keep this reference updated.}
of \cite {AL}) that their critical exponents are also close.

\end{rem}
}

\begin{cor}

Let $F_p: z\mapsto z^2+c_p$ be the straightening of  $f_p$. 
Then
\be
\lim_{p \to \infty} \delta_{\crit}(F_p)=1.
\ee

\end{cor}

\begin{pf}

By Lemma \ref{class}, the germ $f_p$ has a quadratic-like representative $f_p: U_p\ra V_p$
with a big modulus: $\mod(V_p\sm U_p) \to \infty $ as $p\to \infty$. 
Hence for $p$ large, there is a  quasiconformal conjugacy between $f_p$ and $F_p$ with
a small dilatation.
This easily implies (see for instance Lemma 3.15 of \cite {AL}) 
that $f_p$ and $F_p$ have close critical exponents.
\end{pf}

%% file: cheb.tex
\section{Lebesgue measure of the Julia set}

Below $f=f_p :U \to V$ will be the fixed point of nearly Chebyshev
renormalization of period $p$, and $f'=f^p : U'\ra V'$ will be its
pre-renormalization as constructed in Lemma \ref{class}.  Thus,
$f'(z) =\lambda f(\lambda^{-1} z)$, where $\la=\la_p \in (-1,0)$
is the scaling factor of $f$.
We let as above $A=V \setminus U$, $A'=V' \setminus U'$.
Furthermore, we let $U^k=\lambda^k U$, $V^k=\lambda^k V$, and
$A^k=V^k \setminus U^k$.

We will need the following combinatorial lemma.

\begin{lemma} \label {ukvk}

Let
$$
u^k=\sup\, \Xi_\delta(V^k \xleftarrow[U \setminus V^k]{+} A^k),
$$
$$
v^k_j=\sup\, \Xi_\delta^{[j]}(A^k \xleftarrow[U \setminus V^{k+1}]{+} A^k),
$$
$$
v^k=\lim_{j \to \infty} v^k_j=
\sup\, \Xi_\delta(A^k \xleftarrow[U \setminus V^{k+1}]{+} A^k).                   
$$                                                                      
Then
\be \label {71}
u^{k+1} \leq \sup\, \Xi_\delta(V' \xleftarrow[U \setminus V']{+} A')+
u^k (1+v^k) \sup\, \Xi_\delta(A \xleftarrow[U \setminus V']{} A')
(1+\sup\, \Xi_\delta(V' \xleftarrow[U \setminus V']{} U \setminus
V')),
\ee
\be \label {72}
v^k_{j+1} \leq (1+v^k_j)
u^k (1+\sup\, \Xi_\delta(A \xleftarrow[U \setminus V']{} U \setminus V')).
\ee

\end{lemma}

\begin{pf}

Let $B^k=U \setminus (A^k \cup V^{k+1})$.
Let us prove the first estimate.

We can decompose $V^{k+1} \xleftarrow[U \setminus V^{k+1}]{+} A^{k+1}$ into
two groups:
$$
V^{k+1} \xleftarrow[B^k]{+} A^{k+1},
$$
which takes into account the orbits that do not land at the
annulus $A^k$, and
$$
V^{k+1} \xleftarrow[B^k]{} A^k \xleftarrow[U \setminus V^{k+1}]{} A^k
\xleftarrow[B^k]{} A^{k+1},
$$
which accounts for the orbits landing at $A^k$ and 
marks the first and the last landings.  Thus
\begin{align}  \label {11}
u^{k+1}=
&\sup\,\Xi_\delta(V^{k+1} \xleftarrow[U \setminus V^{k+1}]{+} A^{k+1}) \leq
\sup\,\Xi_\delta(V^{k+1} \xleftarrow[B^k]{+} A^{k+1})\\
\nonumber
&+\left (\sup\,\Xi_\delta(V^{k+1} \xleftarrow[B^k]{} A^k)
\cdot \sup\,\Xi_\delta(A^k \xleftarrow[U \setminus V^{k+1}]{} A^k)
\cdot \sup\,\Xi_\delta(A^k \xleftarrow[B^k]{} A^{k+1}) \right ).
\end{align}
Notice that
\be \label {12}
\Xi_\delta(V^{k+1}
\xleftarrow[B^k]{+} A^{k+1})(x)=
\Xi_\delta(V' \xleftarrow[U \setminus V']{+} A')(\lambda^{-k}(x)),
\ee
\be \label {20}
\Xi_\delta(A^k \xleftarrow[B^k]{} A^{k+1})(x)=
\Xi_\delta(A \xleftarrow[U \setminus V']{} A')(\lambda^{-k}x),
\ee
and
\be \label {32}
\sup\,\Xi_\delta(A^k \xleftarrow[U \setminus V^{k+1}]{} A^k)=
1+\sup\,\Xi_\delta(A^k \xleftarrow[U \setminus V^{k+1}]{+} A^k)=1+v^k,
\ee
the $1$ accounting for trivial orbits.
Plugging (\ref {12}) -- (\ref {32}) into (\ref {11}) we get
\be \label {28}
u^{k+1} \leq \sup\, \Xi_\delta(V' \xleftarrow[U \setminus V']{+} A')+
(1+v^k) \sup\, \Xi_\delta(A \xleftarrow[U \setminus V']{} A')
\sup\,\Xi_\delta(V^{k+1} \xleftarrow[B^k]{} A^k).
\ee

We can decompose $V^{k+1} \xleftarrow[B^k]{} A^k$ into two groups,
$$
V^{k+1} \xleftarrow[U \setminus V^k]{} A^k,\quad \text {and} \quad V^{k+1}
\xleftarrow[B^k]{} U^k \setminus V^{k+1} \xleftarrow[U \setminus V^k]{} A^k.
$$
Thus
\begin{align} \label {29}
\sup\, \Xi_\delta(V^{k+1} \xleftarrow[B^k]{} A^k) \leq
&\sup\, \Xi_\delta(V^{k+1} \xleftarrow[U \setminus V^k]{} A^k)\\
\nonumber
&+\sup\, \Xi_\delta(V^{k+1} \xleftarrow[B^k]{} U^k \setminus V^{k+1})
\sup\, \Xi_\delta(U^k \setminus V^{k+1} \xleftarrow[U \setminus V^k]{}
A^k).
\end{align}

Notice that
\be \label {15}
\Xi_\delta(V^{k+1} \xleftarrow[B^k]{} U^k \setminus V^{k+1})(x)=
\Xi_\delta(V' \xleftarrow[U \setminus V']{} U \setminus V')(\lambda^{-k}x),
\ee
\begin{align} \label {16}
\max \left \{
\sup\,\Xi_\delta(V^{k+1} \xleftarrow[U \setminus V^k]{} A^k),\,
\sup\,\Xi_\delta(U^k \setminus V^{k+1} \xleftarrow[U \setminus V^k]{}
A^k),\,
\sup\,\Xi_\delta(A^k \xleftarrow[U \setminus V^k]{+} A^k) \right \}&\\
\nonumber
=\sup\,\Xi_\delta(V^k \xleftarrow[U \setminus V^k]{+} A^k)=u&^k.
\end{align}
Plugging (\ref {15}) and (\ref {16}) into (\ref {29}), and plugging the
resulting expression for $\sup\, \Xi_\delta(V^{k+1} \xleftarrow[B^k]{} A^k)$
into (\ref {28}) gives (\ref {71}).

Let us prove the second estimate.
We will omit the truncation parameter ($j$ or $j+1$).

We can rewrite $A^k \xleftarrow[U \setminus V^{k+1}]{+} A^k$ as
$A^k \xleftarrow[B^k]{+} A^k
\xleftarrow[U \setminus V^{k+1}]{} A^k$.  Thus
\be \label {31}
\sup\,\Xi_\delta(A^k \xleftarrow[U \setminus V^{k+1}]{+} A^k) \leq
\sup\,\Xi_\delta(A^k \xleftarrow[B^k]{+} A^k)
\sup\,\Xi_\delta(A^k \xleftarrow[U \setminus V^{k+1}]{} A^k).
\ee
Plugging (\ref {32}) into
(\ref {31}) we get
\be \label {33}
v^k \leq (1+v^k) \sup\,\Xi_\delta(A^k \xleftarrow[B^k]{+} A^k).
\ee

We can split $A^k \xleftarrow[B^k]{+} A^k$ into
two groups:
$A^k \xleftarrow[U \setminus V^k]{+} A^k$ and
$A^k \xleftarrow[B^k]{} U^k \setminus V^{k+1}
\xleftarrow[U \setminus V^k]{} A^k$.
Thus
\begin{align} \label {35}
\sup\,\Xi_\delta(A^k \xleftarrow[B^k]{+} A^k)
\leq &\sup\,\Xi_\delta(A^k \xleftarrow[U \setminus V^k]{+} A^k)\\
\nonumber
&+\sup\,\Xi_\delta(A^k \xleftarrow[B^k]{} U^k
\setminus V^{k+1})
\sup\,\Xi_\delta(U^k \setminus V^{k+1}
\xleftarrow[U \setminus V^k]{} A^k).
\end{align}

Notice that
\be \label {39}
\Xi_\delta(A^k \xleftarrow[B^k]{} U^k
\setminus V^{k+1})(x)=\Xi_\delta(A \xleftarrow[U \setminus V']{} U
\setminus V')(\lambda^{-k}x),
\ee
Plugging (\ref {39}) and (\ref {16}) into (\ref {35}) we get
\be \label {42}
\sup\,\Xi_\delta(A^k \xleftarrow[B^k]{+} A^k)
\leq u^k+u^k\sup\,\Xi_\delta(A \xleftarrow[U \setminus V']{} U
\setminus V').
\ee
Plugging (\ref {42}) into (\ref {33}) gives (\ref {72}).
\end{pf}

\bigskip

\noindent{\it Proof of Theorem \ref {thmB}.}
By Lemma~\ref {prevlemm}, there exists
$K \equiv K(2)>0$ such that if one takes
$\rho$ sufficiently small, then for all $p$ sufficiently large we have
\be\label{raz}
\sup\, \Xi_2(V' \xleftarrow[U \setminus V']{+} A')<\frac {1} {100},
\ee

\be\label{dva}
\sup\, \Xi_2(A \xleftarrow[U \setminus V']{} A')<\frac {1} {5K+5},
\ee
\be\label{tri}
\sup\, \Xi_2(V' \xleftarrow[U \setminus V']{} U \setminus V')<2 K,
\ee

\be\label{chetyre}
\sup\, \Xi_2(A \xleftarrow[U \setminus V']{} U \setminus V')<\frac {1}
{100}.
\ee

Let us show by induction that for every $k \geq 0$ we have
\be \label {u_k}
u^k \leq \frac {1} {10},
\ee
where $u^k$ is as in Lemma~\ref {ukvk}.
Notice that $u^0=0$, so (\ref {u_k}) holds for $k=0$.  Assuming that (\ref
{u_k}) holds for some $k$, notice that (\ref {72}) and (\ref{chetyre}) imply
$$
v^k_{j+1} \leq \frac {1} {5} (1+v^k_j),
$$
for $j \geq -1$.  Since $v^k_{-1}=0$, this implies by induction
that $v^k_j \leq \frac {1} {4}$ for every $j \geq -1$, so
$v^k=\lim_{j \to \infty} v^k_j \leq \frac {1} {4}$. 
 By (\ref {71}) and (\ref{raz}) -- (\ref{tri}), we have
$$
u^{k+1} \leq \frac {1} {100}+\frac {1} {10} \frac {5} {4} \frac {1} {5K+5}
(2 K+1) \leq \frac {1} {10}.
$$
By induction, (\ref {u_k}) holds for all $k \geq 0$.

Let
$$
X_k=\cup_{r \geq 1} f^{-r} V^k.
$$
Then
$$
\area (X_k \cap A^k)=\int_{A^k} 1_{X^k} dx=\int_{V^k}
\Xi_2(V^k \xleftarrow[U \setminus V^k]{+} A^k) dx \leq
\int_{V^k} u^k dx \leq \frac {1} {10} \area (V^k).
$$
Notice that $X^k \cap V^k=U^k \cup (X^k \cap A^k)$.  Thus,
\be \label {X_k}
\frac {\area (X^k \cap V^k)} {\area V^k} \leq \frac {1} {10}+\frac
{\area U^k } {\area V^k} \leq \frac {1} {5},
\ee
where we have used that $\area U^k \leq \frac {1} {10} \area V^k$,
 which holds since
$\mod (V^k \setminus \overline {U^k})=\mod A $ is big for large $p$ 
(by Lemma \ref{class}).

\medskip

The conclusion of the argument is standard.  Let
$$
X=\{x \in J(f),\, 0 \in \omega(x)\}.
$$
Notice that $X$ is fully invariant: $X=f^{-1}(X)=f(X)$.
By \cite{typical},  for almost every $x \in J(f)$,
$\omega(x) \subset \omega(0)$.  Since $\omega(0)$ is a minimal set
containing $0$, we conclude that $\area X=\area J(f)$.  Let us show
that $\area X=0$.

Assume that this is not the case.  By the Lebesgue Density Points Theorem,\
there exists a density point $x \in X$.  Let $r_k \geq 0$ be minimal
such that $f^{r_k}(x) \in V^k$.  We may assume that $x$ is not a preimage of
$0$, so that $r_k \to \infty$.  Let $W^k$ be the connected component of
$f^{-r_k}(V^k)$ containing $x$.  Then $f^{r_k}:W^k \to V^k$ admits a
univalent extension onto $\V^k \equiv \lambda^k \V$,
and since $\mod (\V^k
\setminus \overline {V^k})$ is big, it has distortion bounded by $2$.  It
also follows that $W^k$ contains a round disk of radius $\frac
{1} {10} \diam (W^k)$.  Since $r_k \to \infty$ and $W^k \subset f^{-r_k}(V)$,
$\limsup W^k \subset K(f)$.  Since $K(f)$ has empty interior, we conclude
that $\diam (W^k) \to 0$.  Notice that
$$
\frac {\area(V^k \setminus X)} {\area V^k } \leq 10 \frac {\area
(W^k \setminus X^k)} {\area W^k } \leq 1000 \frac
{\area (\D_{\diam(W^k)}(x) \setminus X)} {\area (\D_{\diam(W^k)}(x)},
$$
and since $x$ is a density point of $X$, we have
\be \label {X}
               \frac {\area(V^k \setminus X)} {\area (V^k)} \to 0.
\ee
Obviously, $X \subset X_k$, so (\ref {X}) and (\ref {X_k}) give the desired
contradiction.
\qed

%% file: main.bbl
\begin{thebibliography}{BKNS}

\bibitem[AL]{AL} A. Avila and M. Lyubich.  Hausdorff dimension and
conformal measures of Feigenbaum Julia sets. Manuscript 2003. 

\bibitem[B]{B} C.J. Bishop.
Minkowski dimension and the Poincar\'e exponent.
Michigan Math. J. 43 (1996), no. 2, 231--246.

\bibitem[DH]{DH}
 A. Douady and J.H. Hubbard. On the dynamics of polynomial-like maps.
    Ann. Sc. \'{E}c. Norm. Sup., v. 18 (1985), 287-343.

\bibitem[L1]{typical} M. Lyubich. Typical behavior of trajectories of a rational mapping of the sphere. 
Dokl. Akad. Nauk SSSR, v. 268 (1982), 29--32.

\bibitem[L2]{puzzle} M. Lyubich. Dynamics of quadratic polynomials, I-II.
  Acta Math., v. 178 (1997), 185-297.

\bibitem[Mc1]{McM1} C. McMullen.
 Complex dynamics and renormalization.
  Annals of Math. Studies, v. 142, Princeton University Press,  1994.  

\bibitem[Mc2]{McM2} C. McMullen.
 Renormalization and 3-manifolds which fiber over the circle.
Annals of Math. Studies, v. 135, Princeton University Press, 1996.

\bibitem[MS]{MS} W. de Melo and S. van Strien.  One-dimensional dynamics. 
Springer, 1993.

\bibitem[Sh]{Sh} M. Shishikura. 
The Hausdorff dimension of the boundary of the Mandelbrot set and Julia
sets.  Ann. of Math. (2) 147 (1998), no. 2, 225--267.

\bibitem[S1]{S-conformal} D. Sullivan. Conformal dynamical
   systems. Geometric dynamics (Rio de Janeiro, 1981), 725--752, Lecture
   Notes in Math., 1007, Springer, Berlin, 1983.

\bibitem[S2]{S} D. Sullivan.  Bounds, quadratic differentials, and
renormalization conjectures.  AMS Centennial Publications.  {\bf 2}: Mathematics into
Twenty-first Century (1992).  

\bibitem[Y]{Y} B. Yarrington.
 Local connectivity and Lebesgue measure of polynomial Julia sets.
Thesis, Stony Brook 1995.

\end{thebibliography}
